\documentclass[12pt,reqno]{amsart}
\usepackage{amsmath,amssymb,amsfonts,amscd,latexsym,amsthm,mathrsfs,verbatim}
\usepackage[unicode]{hyperref}
\usepackage{fouridx}

\renewcommand{\d}{{\mathrm d}}
\newcommand{\m}{{\mathrm m}}
\renewcommand{\Re}{\operatorname{Re}}
\renewcommand{\Im}{\operatorname{Im}}

\begin{document}

\title[Exercising in complex Mahler measures]{Exercising in complex Mahler measures:\\diamonds are not forever}

\author{Berend Ringeling}
\address{Department of Mathematics, IMAPP, Radboud University, PO Box 9010, 6500~GL Nij\-megen, Netherlands}
\email{b.ringeling@math.ru.nl}

\author{Wadim Zudilin}
\address{Department of Mathematics, IMAPP, Radboud University, PO Box 9010, 6500~GL Nij\-megen, Netherlands}
\email{w.zudilin@math.ru.nl}

\date{29 September 2021}

\begin{abstract}
Recently, Hang Liu and Hourong Qin came up with a numerical observation about the relation between the Mahler measures of one hyperelliptic and two elliptic families.
The discoverers foresee a proof of the identities ``by extending ideas in'' two papers of M.~Lal\'\i n and Gang Wu, the ideas based on a theorem of S.~Bloch and explicit diamond-operation calculations on the underlying curves.
We prove the relation using the already available diamond-free methodology.
While finding such relations for the Mahler measures remains an art, proving them afterwards is mere complex (analysis) exercising.
\end{abstract}


\thanks{The work is supported by NWO grant OCENW.KLEIN.006.}

\maketitle

It is still a mystery what the (logarithmic) Mahler measure
\begin{align*}
\m(P(x_1,\dots,x_k))
&=\frac1{(2\pi i)^k}\idotsint_{|x_1|=\dots=|x_k|=1}\log|P(x_1,\dots,x_k)|\,\frac{\d x_1}{x_1}\dotsb\frac{\d x_k}{x_k}
\end{align*}
of a nonzero $k$-variable (Laurent) polynomial $P(x_1,\dots,x_k)\in\mathbb Z[x_1,\dots,x_k]$ is when $k\ge2$.
The single-variable case is quite special, because Jensen's classical formula expresses $\m(P(x))$ as a sum of the logarithms of algebraic integers;
however we do not yet know whether it is bounded away from 0 when $\m(P(x))>0$.
A significant progress was done in the 1980s and 1990s towards understanding the nature of two-variable Mahler measures $\m(P(x,y))$ at least for the cases when $P(x,y)=0$ represents a curve of small genus.
In his work \cite{Bo98}, D.~Boyd tabulated examples of families of such Mahler measures which are ultimately linked to the $L$-values at 2 of the algebraic varieties related to $P(x,y)=0$.
In several cases he observed explicit $\mathbb Z$-linear relations between the Mahler measures from two different families;
for example, he noticed that
\begin{equation}
\m\big(Q_k(X,Y)\big)=\begin{cases}
2\m\big(P_{k-4}(x,y)\big) &\text{for}\; 0\le k\le 4, \\
\phantom1\m\big(P_{k-4}(x,y)\big) &\text{for}\; k\le-1,
\end{cases}
\label{boyd}
\end{equation}
where
\begin{align*}
Q_k(X,Y)
&=Y^2+(X^4+kX^3+2kX^2+kX+1)Y+X^4,
\\
P_\lambda(x,y)
&=(x+1)y^2+(x^2-(\lambda+2)x+1)y+x(x+1).
\end{align*}
Here $Q_k(X,Y)=0$ is generically a hyperelliptic curve, while $P_\lambda(x,y)=0$ is elliptic.
The expectation \eqref{boyd} was settled by M.\,J.~Bertin and one of the present authors in \cite{BZ16} who also gave other examples of relations between hyperelliptic and elliptic Mahler measures in~\cite{BZ17}.
Later in \cite{LW19} M.~Lal\'\i n and G.~Wu gave a different proof of \eqref{boyd} based on a theorem of S.~Bloch and explicit diamond-operation calculations on the underlying (hyper)elliptic curves.
More recently in \cite{LW20} they used the method to prove some newer numerical observations of this type.
At the same time H.~Liu and H.~Qin made a systematic search of such relations in \cite{LQ21} with the following nice example singled out from their collection:
\begin{equation}
\m\big(Q_{\lambda+4}(X-1,Y)\big)
=\begin{cases}
\m(R_\lambda(x,y)) &\text{if}\; \lambda\le-5, \\
\frac12\big(\m(R_\lambda(x,y))+\m(P_\lambda(x,y))\big) &\text{if}\; \lambda\ge13, \\
\end{cases}
\label{main}
\end{equation}
where $Q_k(X,Y)$, $P_\lambda(x,y)$ are as above, while
\[
R_\lambda(x,y)
=x+x^{-1}+y+y^{-1}+\lambda
\]
is the `most featured' family in the Mahler measure literature.
The relations in \eqref{main} not only involve a bunch of families, all on the original Boyd's list in \cite{Bo98}, but also mixes up \emph{three} different Mahler measures.
Interestingly enough, the discoverers of \eqref{main} foresee a proof of the relations ``by extending ideas in'' \cite{LW19,LW20} rather than through the diamond-free machinery in \cite{BZ16,BZ17}.
The goal of this note is to demonstrate that the latter methodology applies well to~\eqref{main}, without new inventions.
While finding relations like \eqref{boyd} and \eqref{main} for the Mahler measures remains an art,
proving them afterwards is not more than complex (analysis) exercising.

\medskip
Observe that
\[
Q_{\lambda+4}(X-1,Y)
=X^8\cdot\big(y^2+(2x^2+\lambda x+1)y+x^4\big)\Big|_{x=(X-1)/X^2,\,y=Y/X^4},
\]
hence
\begin{align}
q(\lambda)
&=\m\big(Q_{\lambda+4}(X-1,Y)\big)
\nonumber\\
&=\frac1{(2\pi i)^2}\int_{|y|=1}\frac{\d y}y\int_{|X|=1}\log|y^2+(2x^2+\lambda x+1)y+x^4|\,\Big|_{x=(X-1)/X^2}\,\frac{\d X}X.
\label{m1}
\end{align}
The zeros of the polynomial $y^2+(2x^2+\lambda x+1)y+x^4$ are
\[
y_{\pm}(x)=(\lambda x+1)\bigg(\frac12+\frac{x^2}{\lambda x+1}\pm\sqrt{\frac14+\frac{x^2}{\lambda x+1}}\bigg).
\]
If we use $X=e^{-2\pi it}$ where $t\in[-1/2,1/2]$ as parameterisation of the circle $|X|=1$, so that $x(t)=(X-1)/X^2=e^{2\pi it}(1-e^{2\pi it})$, then one checks that in the case $\lambda\ge13$,
\[
|y_-(x(t))|\le|y_-(x(\pm1/2))|\le1 \quad\text{and}\quad |y_+(x(t))|\ge|y_+(x(0))|=1
\]
for all $t\in[-1/2,1/2]$; the behaviour of $|y_\pm(x)|$ in the case $\lambda\le-4$ is somewhat subtler but we also have
\[
|y_-(x(t))|\le1 \quad\text{and}\quad |y_+(x(t))|\ge1
\]
for all $t\in[-1/2,1/2]$.
Therefore, applying to the external integral in \eqref{m1} Jensen's formula we deduce, for both $\lambda\ge13$ and $\lambda\le-4$, that
\begin{align*}
q(\lambda)
&=\int_{-1/2}^{1/2}\log|y_+(x(t))|\,\d t
=2\int_0^{1/2}\log|y_+(x(t))|\,\d t
\displaybreak[2]\\
&=2\Re\int_0^{1/2}\log\bigg((\lambda x+1)\bigg(\frac12+\frac{x^2}{\lambda x+1}+\sqrt{\frac14+\frac{x^2}{\lambda x+1}}\bigg)\bigg)\bigg|_{x=x(t)}\d t.
\end{align*}
Then for $\lambda>13$ and $\lambda<-4$ we obtain
\begin{align}
\frac{\d q(\lambda)}{\d\lambda}
&=2\Re\int_0^{1/2}
\frac{x}{(1+\lambda x)\sqrt{1+4x^2/(1+\lambda x)}}\bigg|_{x=e^{2\pi it}(1-e^{2\pi it})}\,\d t
\nonumber\\
&=\Re\frac1{\pi i}\int_{\substack{|z|=1\\ \Im z>0}}
\frac{x}{(1+\lambda x)\sqrt{1+4x^2/(1+\lambda x)}}\bigg|_{x=z(1-z)}\,\frac{\d z}z.
\label{der_q}
\end{align}

In order to flatten the resulting arc integral \eqref{der_q}, let us analyse the location of singularities of the integrand.
Since $x=z(1-z)$ maps the real interval $-1<z<1$ onto $-2<x<1/4$, we first need to look after the zeros of the cubic polynomial $(1+\lambda x)(1+\lambda x+4x^2)$, explicitly
$$
x_0=x_0(\lambda)=-\frac1\lambda, \;\;
x_1=x_1(\lambda)=-\frac{\lambda+\sqrt{\lambda^2-16}}8
\;\;\text{and}\;\;
x_2=x_2(\lambda)=-\frac{\lambda-\sqrt{\lambda^2-16}}8,
$$
on the latter interval.
Furthermore, notice that the interval $-1<z<1/2$ is bijectively mapped onto $-2<x<1/4$ with the inverse mapping $z=(1-\sqrt{1-4x})/2$, while the interval $1/2<z<1$ is bijectively mapped onto $0<x<1/4$ reversed(!) with the inverse mapping $z=(1+\sqrt{1-4x})/2$.

When $\lambda<-5$, we have $0<x_0<x_1<1/4$ and $x_2>1$.
The corresponding $z$-singularities
\[
z_1=\frac{1-\sqrt{1-4x_0}}2, \;\;
z_2=\frac{1-\sqrt{1-4x_1}}2, \;\;
z_3=\frac{1+\sqrt{1-4x_1}}2 \;\;\text{and}\;\;
z_4=\frac{1+\sqrt{1-4x_0}}2
\]
(in this order!) are all located on $0<z<1$; analysing the sign of the square root in the denominator in~\eqref{der_q} on each of the $z$-intervals $(-1,z_1)$, $(z_1,z_2)$, $(z_2,z_3)$, $(z_3,z_4)$ and $(z_4,1)$ we find out that
\begin{align*}
&
\int_{\substack{|z|=1\\ \Im z>0}}
\frac{x}{(1+\lambda x)\sqrt{1+4x^2/(1+\lambda x)}}\bigg|_{x=z(1-z)}\,\frac{\d z}z
\\ &\quad
=-\bigg(\int_{-1}^{z_1}+\int_{z_1}^{z_2}+\int_{z_2}^{z_3}-\int_{z_3}^{z_4}+\int_{z_4}^1\bigg)
\frac{x}{(1+\lambda x)\sqrt{1+4x^2/(1+\lambda x)}}\bigg|_{x=z(1-z)}\,\frac{\d z}z.
\end{align*}
Thus, the integral \eqref{der_q} transforms for $\lambda<-5$ into
\begin{align}
\frac{\d q(\lambda)}{\d\lambda}
&=-\frac1{\pi i}\bigg(\int_{z_1}^{z_2}-\int_{z_3}^{z_4}\bigg)
\frac{x}{(1+\lambda x)\sqrt{1+4x^2/(1+\lambda x)}}\bigg|_{x=z(1-z)}\,\frac{\d z}z
\nonumber\displaybreak[2]\\
&=-\frac1{2\pi i}\bigg(\int_{x_0}^{x_1}\frac{(1+\sqrt{1-4x})\,\d x}{(1+\lambda x)\sqrt{1+4x^2/(1+\lambda x)}\sqrt{1-4x}}
\nonumber\\ &\qquad
-\int_{x_0}^{x_1}
\frac{(1-\sqrt{1-4x})\,\d x}{(1+\lambda x)\sqrt{1+4x^2/(1+\lambda x)}\sqrt{1-4x}}\bigg)
\nonumber\displaybreak[2]\\
&=-\frac1\pi\int_{-1/\lambda}^{-(\lambda+\sqrt{\lambda^2-16})/8}
\frac{\d x}{\sqrt{-(1+\lambda x)(1+\lambda x+4x^2)}}.
\label{der-}
\end{align}

When $\lambda>13$, we get $x_1<-2<x_2<x_0<0$, so that the $z$-singularities on the interval $-1<z<1$ are
\[
z_1=\frac{1-\sqrt{1-4x_2}}2 \;\;\text{and}\;\;
z_2=\frac{1-\sqrt{1-4x_0}}2.
\]
The analysis as above brings us to
\begin{align*}
&
\int_{\substack{|z|=1\\ \Im z>0}}
\frac{x}{(1+\lambda x)\sqrt{1+4x^2/(1+\lambda x)}}\bigg|_{x=z(1-z)}\,\frac{\d z}z
\\ &\quad
=-\bigg(\int_{-1}^{z_1}-\int_{z_1}^{z_2}+\int_{z_2}^1\bigg)
\frac{x}{(1+\lambda x)\sqrt{1+4x^2/(1+\lambda x)}}\bigg|_{x=z(1-z)}\,\frac{\d z}z
\end{align*}
and we conclude that, for $\lambda>13$,
\begin{align}
\frac{\d q(\lambda)}{\d\lambda}
&=\frac1{\pi i}\int_{z_1}^{z_2}\frac{x}{(1+\lambda x)\sqrt{1+4x^2/(1+\lambda x)}}\bigg|_{x=z(1-z)}\,\frac{\d z}z
\nonumber\displaybreak[2]\\
&=\frac1{2\pi i}\int_{x_2}^{x_0}\frac{(1+\sqrt{1-4x})\,\d x}{(1+\lambda x)\sqrt{1+4x^2/(1+\lambda x)}\sqrt{1-4x}}
\nonumber\displaybreak[2]\\
&=\frac1{2\pi}\bigg(\int_{-(\lambda-\sqrt{\lambda^2-16})/8}^{-1/\lambda}\frac{\d x}{\sqrt{-(1+\lambda x)(1+\lambda x+4x^2)}}
\nonumber\\ &\qquad
+\int_{-(\lambda-\sqrt{\lambda^2-16})/8}^{-1/\lambda}\frac{\d x}{\sqrt{-(1+\lambda x)(1+\lambda x+4x^2)(1-4x)}}\bigg).
\label{der+}
\end{align}

On the other hand, consider $p(\lambda)=\m(P_\lambda(x,y))$ and $r(\lambda)=\m(R_\lambda(x,y))=r(|\lambda|)$.
When $\lambda>5$, the substitution
\[
y =\frac{-\lambda - \sqrt{\lambda^2-16} + (8 - \lambda^2 - \lambda \sqrt{\lambda^2-16})x}
{-\lambda + \sqrt{\lambda^2-16} + (8 - \lambda^2 + \lambda\sqrt{\lambda^2-16})x}
\]
and the internal parameterisation $\lambda=2(1+\mu^2)/\mu$ where $0<\mu<1/2$ leads us to
\begin{align}
&
\int_{-(\lambda-\sqrt{\lambda^2-16})/8}^{-1/\lambda}
\frac{\d x}{\sqrt{-(1+\lambda x)(1+\lambda x+4x^2)(1-4x)}}
\nonumber\\  &\quad
= \frac{\mu}{2} \int_0^1 \frac{\d y}{\sqrt{y(1-y)(1+2\mu-\mu^3(2+\mu)y)}}
= \frac{\pi\mu}{2 \sqrt{1+2\mu}}\,{}_2F_1\bigg(\begin{matrix} \frac12, \, \frac12 \\ 1 \end{matrix} \biggm|\frac{\mu^3(2+\mu)}{1+2\mu} \bigg)
\nonumber\\  &\quad
= \frac{\pi\mu}{2(1+4\mu+\mu^2)}\,{}_2 F_1\bigg(\begin{matrix} \frac13, \, \frac23 \\ 1 \end{matrix} \biggm|\frac{27\mu (1+\mu)^4}{2(1+4\mu+\mu^2)^3} \bigg)
= \frac\pi{\lambda+8}\,{}_2F_1\bigg(\begin{matrix} \frac13, \, \frac23 \\ 1 \end{matrix} \biggm|\frac{27(\lambda+4)^2}{(\lambda+8)^3}\bigg)
\nonumber\\  \intertext{(we apply \cite[eq.~(5)]{BZ16})} 
&\quad
= \pi \,\frac{\d p(\lambda)}{\d \lambda}.
\label{J1}
\end{align}
For the integral in \eqref{der-}, when $\lambda<-5$, we substitute
\[
z = \frac{x +(\lambda + \sqrt{\lambda^2-16})/8}{-1/\lambda +(\lambda + \sqrt{\lambda^2-16})/8}
\]
and use the parameterisation $\lambda=2(1+\mu^2)/\mu$ now with $-1/2<\mu<0$:
\begin{align}
&
\int_{-1/\lambda}^{-(\lambda+\sqrt{\lambda^2-16})/8}\frac{\d x}{\sqrt{-(1+\lambda x)(1+\lambda x+4x^2)}}
\nonumber\\  &\quad
= -\frac{\mu}{2\sqrt{1-\mu^4}} \int_0^1 \frac{\d z}{\sqrt{z(1-z)(1+z\mu^4/(1-\mu^4))}}
=  -\frac{\pi\mu}{2\sqrt{1-\mu^4}} \,{}_2F_1\bigg(\begin{matrix} \frac12, \, \frac12 \\ 1 \end{matrix} \biggm|\frac{-\mu^4}{1-\mu^4} \bigg)
\nonumber\\  &\quad
= -\frac{\pi\mu}{2(1+\mu^2)}\,{}_2F_1\bigg(\begin{matrix} \frac12, \, \frac12 \\ 1 \end{matrix} \biggm|\frac{4\mu^2}{(1+\mu^2)^2} \bigg)
= - \int_0^1 \frac{\d t}{\sqrt{t(1-t)(\lambda^2-16t)}}
\nonumber\\  \intertext{(see, for example, \cite[Section 3.5]{BZ20})} &\quad
=\pi\,\frac{\d r(\lambda)}{\d\lambda};
\label{J2}
\end{align}
a similar manipulation results in
\begin{equation}
\int_{-(\lambda-\sqrt{\lambda^2-16})/8}^{-1/\lambda}
\frac{\d x}{\sqrt{-(1+\lambda x)(1+\lambda x+4x^2)}}
=\int_0^1\frac{\d t}{\sqrt{t(1-t)(\lambda^2-16t)}}
=\pi\,\frac{\d r(\lambda)}{\d\lambda}
\label{J3}
\end{equation}
when $\lambda>5$.
Proofs of intermediate hypergeometric transformations
\[
\frac1{\sqrt{1+2\mu}}\,{}_2F_1\bigg(\begin{matrix} \frac12, \, \frac12 \\ 1 \end{matrix} \biggm|\frac{\mu^3(2+\mu)}{1+2\mu} \bigg)
= \frac1{1+4\mu+\mu^2}\,{}_2 F_1\bigg(\begin{matrix} \frac13, \, \frac23 \\ 1 \end{matrix} \biggm|\frac{27\mu (1+\mu)^4}{2(1+4\mu+\mu^2)^3} \bigg)
\]
and
\[
\frac1{\sqrt{1-\mu^4}} \,{}_2F_1\bigg(\begin{matrix} \frac12, \, \frac12 \\ 1 \end{matrix} \biggm|\frac{-\mu^4}{1-\mu^4} \bigg)
=\frac1{1+\mu^2}\,{}_2F_1\bigg(\begin{matrix} \frac12, \, \frac12 \\ 1 \end{matrix} \biggm|\frac{4\mu^2}{(1+\mu^2)^2} \bigg)
\]
are straightforward: one can use a computer algebra system to the check the coincidence of linear differential equations for both sides or consult alternatively with~\cite{CZ19}.

Comparing \eqref{der-}, \eqref{der+} with \eqref{J1}--\eqref{J3} we conclude that
\begin{alignat*}{2}
\frac{\d q(\lambda)}{\d\lambda}
&=\frac{\d r(\lambda)}{\d\lambda}
&\quad\text{if}\; \lambda&<-5,
\\
\frac{\d q(\lambda)}{\d\lambda}
&=\frac12\bigg(\frac{\d r(\lambda)}{\d\lambda}+\frac{\d p(\lambda)}{\d\lambda}\bigg)
&\quad\text{if}\; \lambda&>13.
\end{alignat*}
Taking the anti-derivatives, using the continuity of all $q(\lambda),r(\lambda),p(\lambda)$ on the closed intervals $\lambda\le-5$, $\lambda\ge13$ and the asymptotics $\log|\lambda|+O(1)$ as $\lambda\to\infty$ for each of the Mahler measures in question, we finally arrive at \eqref{main} wanted.

\end{document}